\documentclass[12pt]{article}

\usepackage{amsmath}

\usepackage[english]{babel}
\usepackage{amsmath}
\usepackage{amsthm,amssymb,latexsym}
\usepackage{epic, eepic, graphics}

\newtheorem{prop}{Proposition}

\newtheorem{cor}[prop]{Corollary}
\newtheorem{theorem}[prop]{Theorem}
\theoremstyle{definition}

\newtheorem{problem}{Problem}

\newcommand\p{\circle*{0.3}}

\title {On uniquely $k$-determined permutations}
\author {Sergey Avgustinovich\footnote{Sobolev Institute of Mathematics, Acad. Koptyug prospect 4,
Novosibirsk 630090, Russia, avgust@math.nsc.ru}\ \ and Sergey
Kitaev\footnote{Institute of Mathematics, Reykjav\'{\i}k University,
Ofanleiti 2, 101 Reykjav\'{\i}k, Iceland, sergey@ru.is}}
\begin{document}
\maketitle
\begin{abstract}
There are several approaches to study occurrences of consecutive
patterns in permutations such as the inclusion-exclusion method, the
tree representations of permutations, the spectral approach and
others. We propose yet another approach to study occurrences of
consecutive patterns in permutations. The approach is based on
considering the graph of patterns overlaps, which is a certain
subgraph of the de Bruijn graph.

While applying our approach, the notion of a uniquely $k$-determined
permutation appears. We give two criteria for a permutation to be
uniquely $k$-determined: one in terms of the distance between two
consecutive elements in a permutation, and the other one in terms of
directed hamiltonian paths in the certain graphs called
path-schemes. Moreover, we describe a finite set of prohibitions
that gives the set of uniquely $k$-determined permutations. Those
prohibitions make applying the transfer matrix method possible for
determining the number of uniquely $k$-determined permutations.
\end{abstract}

\section{Introduction}

A {\em pattern} $\tau$ is a permutation on $\{1,2,\ldots,k\}$. An
occurrence of a {\em consecutive pattern} $\tau$ in a permutation
$\pi=\pi_1\pi_2\ldots\pi_n$ is a word
$\pi_i\pi_{i+1}\ldots\pi_{i+k-1}$ that is {\em order-isomorphic} to
$\tau$. For example, the permutation 253164 contains two occurrences
of the pattern 132, namely 253 and 164. In this paper we deal only
with consecutive patterns, which courses omitting the word
``consecutive'' in defining a pattern to shorten the notation.

There are several approaches in the literature to study the {\em
distribution} and, in particular, {\em avoidance}, of consecutive
patterns in permutations. For example, direct combinatorial
considerations are used in~\cite{Kit0}; the {\em method of
inclusion-exclusion} is used in~\cite{GoulJack,Kit1}; the {\em tree
representations of permutations} are used in~\cite{ElizNoy}; the
{\em spectral theory of integral operators} on $L^{2}([0,1]^{k})$ is
used in~\cite{EhrKitPer}. In this paper we introduce yet another
approach to study occurrences of consecutive patterns in
permutations. The approach is based on considering the {\em graph of
patterns overlaps} defined below, which is a subgraph of the {\em de
Bruijn graph} studied broadly in the literature mainly in connection
with combinatorics on words and graph theory.

Suppose we are interested in the number of occurrences of a pattern
$\tau$ of length $k$ in a permutation $\pi$ of length $n$. To find
this number, we scan $\pi$ from left to right with a ``window'' of
length $k$, that is, we consider
$P_i=\pi_i\pi_{i+1}\ldots\pi_{i+k-1}$ for $i=1,2,\ldots,n-k+1$: if
we meet an occurrence of $\tau$, we register it. Each $P_i$ forms a
pattern of length $k$, and the procedure of scanning $\pi$ gives us
a path in the {\em graph $\mathcal{P}_k$ of patterns overlaps of
order $k$} defined as follows (graphs of patterns/permutations
overlaps appear in \cite{BurKit,Chung,Hur}). The nodes of
$\mathcal{P}_k$ are all $k!$ $k$-permutations, and there is an arc
from a node $a_1a_2\ldots a_k$ to a node $b_1b_2\ldots b_k$ if and
only if $a_2a_3\ldots a_k$ and $b_1b_2\ldots b_{k-1}$ form the same
pattern. Thus, for any $n$-permutation there is a path in
$\mathcal{P}_k$ of length $n-k+1$ corresponding to it. For example,
if $k=3$ then to the permutation $13542$ there corresponds the path
$123\rightarrow 132\rightarrow 321$ in $\mathcal{P}_3$.

Our approach to study the distribution of a consecutive pattern
$\tau$ of length $k$ among $n$-permutations is to take
$\mathcal{P}_k$ and to consider all paths of length $n-k+1$ passing
through the node $\tau$ exactly $\ell$ times, where
$\ell=0,1,\ldots, n-k+1$. Then we could count the permutations
corresponding to the paths. Similarly, for the ``avoidance
problems'' that attracted much attention in the literature, we
proceed as follows: given a set of patterns of length $k$ to avoid,
we remove the corresponding nodes with the corresponding arcs from
$\mathcal{P}_k$, consider all the paths of certain length in the
graph obtained, and then count the permutations of interest.

However, a complication with the approach is that a permutation does
not need to be reconstructible uniquely from the path corresponding
to it. For example, the permutation $13542$ above has the same path
in $\mathcal{P}_3$ corresponding to it as the permutations $23541$
and $12543$. Thus, different paths in $\mathcal{P}_k$ may have
different contributions to the number of permutations with required
properties; in particular, some of the paths in $\mathcal{P}_k$ give
exactly one permutation corresponding to them. We call such
permutations {\em uniquely $k$-determined}. Study of such
permutations is the main concern of the paper, and it should be
considered as the first step in understanding how to use our
approach to the problems described. Also, in our considerations we
assume that all the nodes in $\mathcal{P}_k$ are allowed while
dealing with uniquely $k$-determined permutations, that is, we do
not prohibit any pattern.

The paper is organized as follows. In Section~\ref{sec2} we study
the set of uniquely $k$-determined permutations. In particular, we
give two criteria for a permutation to be uniquely $k$-determined:
one in terms of the distance between two {\em consecutive elements}
in a permutation, and the other one in terms of directed hamiltonian
paths in the certain graphs called {\em path-schemes}. We use the
second criteria to establish (rough) upper and lower bounds for the
number of uniquely $k$-determined permutations. Moreover, given an
integer $k$, we describe a finite set of prohibitions that
determines the set of uniquely $k$-determined permutations. Those
prohibitions make applying the {\em transfer matrix
method}~\cite[Thm. 4.7.2]{Stan} possible for determining the number
of uniquely $k$-determined permutations and we discuss this in
Subsection~\ref{proh}. As a corollary of using the method, we get
that the generating function for the number of uniquely
$k$-determined permutations is rational. Besides, we show that there
are no {\em crucial permutations} in the set of uniquely
$k$-determined permutations. (Crucial objects, in the sense defined
below, are natural to study in infinite sets of objects defined by
prohibitions; for instance, see~\cite{EvdKit} for some results in
this direction related to words.) We consider in more details the
case $k=3$ in Subsection~\ref{k=3}. Finally, in Section~\ref{sec3},
we state several open problems for further research.

\section{Uniquely $k$-determined permutations}\label{sec2}

\subsection{Distance between consecutive elements; a criterion on unique $k$-determinability}

Suppose $\pi=\pi_1\pi_2\ldots\pi_n$ is a permutation and $i<j$. The
{\em distance} $d_{\pi}(\pi_i,\pi_j)=d_{\pi}(\pi_j,\pi_i)$ between
the elements $\pi_i$ and $\pi_j$ is $j-i$. For example,
$d_{253164}(3,6)=d_{253164}(6,3)=2$.

\begin{theorem}\label{cr01}{\rm[}First criterion on unique $k$-determinability{\rm]}
An $n$-permutation $\pi$ is uniquely $k$-determined if and only if
for each $1\leq x< n$, the distance $d_{\pi}(x,x+1)\leq k-1$.
\end{theorem}

\begin{proof} Suppose for an $n$-permutation $\pi$, $d(x,x+1)\geq k$ for some $1\leq x < n$. This
means that $x$ and $x+1$ will never be inside a ``window'' of length
$k$ while scanning consecutive elements of $\pi$. Thus, these
elements are incomparable in $\pi$ in the sense that switching $x$
and $x+1$ in $\pi$ will lead to another permutation $\pi'$ having
the same path in $\mathcal{P}_k$ as $\pi$ has. So, $\pi$ is not
uniquely $k$-determined.

On the other hand, if for each $1\leq x< n$, the distance
$d_{\pi}(x,x+1)\leq k-1$, then the positions of the elements
$1,2,\ldots, n$ are uniquely determined (first we note that the
position of 1 is uniquely determined, then we determine the position
of 2 which is a 1's neighbor in a ``window'' of length $k$, then the
position of 3, etc.) leading to the fact that $\pi$ is uniquely
$k$-determined.
\end{proof}

The following corollary to Theorem~\ref{cr01} is straightforward.

\begin{cor}\label{cr02} An $n$-permutation $\pi$ is not uniquely $k$-determined
if and only if there exists $x$, $1\leq x< n$, such that
$d_{\pi}(x,x+1)\geq k$. \end{cor}

So, to determine if a given $n$-permutation is uniquely
$k$-determined, all we need to do is to check the distance for $n-1$
pairs of numbers: $(1,2)$, $(2,3)$,..., $(n-1,n)$. Also, the
language of uniquely determined $k$-permutations is {\em factorial}
in the sense that if $\pi_1\pi_2\ldots\pi_n$ is uniquely
$k$-determined, then so is the pattern of $\pi_i\pi_{i+1}\ldots
\pi_j$ for any $i\leq j$ (this is a simple corollary to
Theorem~\ref{cr01}).

Coming back to the permutation $13542$ above and using
Corollary~\ref{cr02}, we see why this permutation is not uniquely
$3$-determined ($k=3$): the distance $d_{13542}(2,3)=3=k$.

\subsection{Directed hamiltonian paths in path-schemes; another criterion on unique
$k$-determinability}\label{dirpaths}

\begin{figure}[h]
\begin{center}
\begin{picture}(6,2)
\put(0,0){\put(0,0){\p} \put(2,0){\p} \put(4,0){\p} \put(6,0){\p}
\put(8,0){\p} \put(10,0){\p}

\put(0,0.5){1} \put(2,0.5){2} \put(4,0.5){3} \put(5.5,0.5){4}
\put(7.5,0.5){5} \put(9.6,0.5){6}

\spline(0,0)(2,-1)(4,0)
\spline(2,0)(4,-1)(6,0)\spline(4,0)(6,-1)(8,0)\spline(6,0)(8,-1)(10,0)

\spline(0,0)(4,2.5)(8,0)\spline(2,0)(6,2.5)(10,0) }
\end{picture}
\caption{The path-scheme $P(6,\{2,4\})$.} \label{il}
\end{center}
\end{figure}
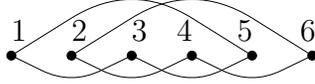

Let $V=\{1,2,\ldots, n\}$ and $M$ be a subset of $V$. A {\em
path-scheme} $P(n,M)$ is a graph $G=(V,E)$, where the edge set $E$
is $\{(x,y)\ |\ |x-y| \in M \}$. See Figure~\ref{il} for an example
of a path-scheme.

Path-schemes appeared in the literature, for example, in connection
with counting independent sets (see~\cite{Kit2}). However, we will
be interested in path-schemes having $M=\{1,2,\ldots, k-1\}$ for
some $k$ (the number of independent sets for such $M$ in case of $n$
nodes is given by the $(n+k)$-th $k$-{\em generalized Fibonacci
number}). Let $\mathcal{G}_{k,n}=P(n,\{1,2,\ldots, k-1\})$, where
$k\leq n$. Clearly, $\mathcal{G}_{k,n}$ is a subgraph of
$\mathcal{G}_{n,n}$.

Any permutation $\pi=\pi_1\pi_2\ldots\pi_n$ determines uniquely a
directed hamiltonian path in $\mathcal{G}_{n,n}$ starting with
$\pi_1$, then going to $\pi_2$, then to $\pi_3$ and so on. The
reverse is also true: given a directed hamiltonian path in
$\mathcal{G}_{n,n}$ we can easily construct the permutation
corresponding to it.

\begin{theorem}\label{cr03}{\rm[}Second criterion on unique $k$-determinability{\rm]}
Let $\Phi$ be a map that sends a uniquely $k$-determined
$n$-permutation $\pi$ to the directed hamiltonian path in
$\mathcal{G}_{n,n}$ corresponding to $\pi^{-1}$. $\Phi$ is a
bijection between the set of all uniquely $k$-determined
$n$-permutations and the set of all directed hamiltonian paths in
$\mathcal{G}_{k,n}$.
\end{theorem}

\begin{proof} Let $\pi$ be a uniquely $k$-determined $n$-permutation.
We claim that the directed hamiltonian path in $\mathcal{G}_{n,n}$
corresponding to $\pi^{-1}$ is actually a directed hamiltonian path
in $\mathcal{G}_{k,n}$. Indeed, suppose the elements $x$ and $x+1$,
$1\leq x <n$, are located in $\pi$ in positions $i$ and $j$
respectively. According to Theorem~\ref{cr01}, $|j-i|\leq k-1$. Now,
$ij$ is a factor in $\pi^{-1}$, and the directed hamiltonian path
corresponding to $\pi^{-1}$ contains the arc from $i$ to $j$, which
is an arc in $\mathcal{G}_{k,n}$. Obviously, $\Phi$ is injective.
Also, it is easy to see how to find the inverse to $\Phi$ mapping a
directed hamiltonian path in $\mathcal{G}_{k,n}$ to a permutation
that, due to Theorem~\ref{cr01}, is uniquely $k$-determined.
\end{proof}

Theorem~\ref{cr03} suggests a quick checking of whether an
$n$-permutation $\pi$ is uniquely $k$-determined or not. One simply
needs to consider $n-1$ differences of the adjacent elements in
$\pi^{-1}$ and check whether at least one of those differences
exceeds~$k-1$ or not. Moreover, one can find the number of uniquely
$k$-determined $n$-permutations by listing them and checking for
each of them the differences of consecutive elements in the manner
described above. Using this approach, one can run a computer program
to get the number of uniquely $k$-determined $n$-permutations for
initial values of $k$ and $n$, which we record in Table~\ref{data}.

\begin{table}[ht]
\begin{center}
   \begin{tabular}{|l|l|}
    \hline
$k=2$ & $1,\ 2,\ 2,\ 2,\ 2,\ 2,\ 2,\ 2,\
2,\ldots$\\
   \hline
$k=3$ & $1,\ 2,\ 6,\ 12,\ 20,\ 34,\ 56,\ 88,\
136,\ldots$\\
   \hline
$k=4$ & $1,\ 2,\ 6,\ 24,\ 72,\ 180,\ 428,\ 1042,\
2512,\ldots$ \\
   \hline
$k=5$ & $1,\ 2,\ 6,\ 24,\ 120,\ 480,\ 1632,\ 5124,\
15860,\ldots$\\
   \hline
$k=6$ & $1,\ 2,\ 6,\ 24,\ 120,\ 720,\ 3600,\ 15600,\
61872,\ldots$\\
   \hline
$k=7$ & $1,\ 2,\ 6,\ 24,\ 120,\ 720,\ 5040,\ 30240,\
159840,\ldots$\\
   \hline
$k=8$ & $1,\ 2,\ 6,\ 24,\ 120,\ 720,\ 5040,\ 40320,\
282240,\ldots$\\
   \hline
   \end{tabular}
\smallskip \caption{The initial values  for the number of uniquely $k$-determined
$n$-permutations.} \label{data}
\end{center}
\end{table}

It is remarkable that the sequence corresponding to the case $k=3$
in Table~\ref{data} appears in~\cite[A003274]{Sloane}, where we
learn that the inverses to the uniquely 3-determined permutations
are called the {\em key permutations} and they appear
in~\cite{Page}. Another sequence appearing in Table~\ref{data} is
\cite[A003274]{Sloane}: 0, 2, 12, 72, 480, 3600, .... In our case,
this is the number of uniquely $n$-determined $(n+1)$-permutations,
$n\geq 1$; in \cite{Sloane}, this is the number of
$(n+1)$-permutations that have 2 predetermined elements non-adjacent
(e.g., for $n=2$, the permutations with say 1 and 2 non-adjacent are
132 and 231). It is clear that both of the last objects are counted
by $n!(n-1)$. Indeed, to create a uniquely $n$-determined
$(n+1)$-permutation, we take any permutation (there are $n!$
choices) and extend it to the right by one element making sure that
the extension is not adjacent to the leftmost element of the
permutation (there are $n-1$ possibilities; here we use
Theorem~\ref{cr01}). On the other hand, to create a ``good''
permutation appearing in \cite{Sloane}, we take any of $n!$
permutations, and insert one of the predetermined elements into any
position not adjacent to the other predetermined element (there are
$(n-1)$ choices). A bijection between the sets of permutations above
is given by the following: Suppose $a$ and $b$ are the predetermined
elements in $\pi=\pi_1\ldots\pi_n$, and $\pi_i=a$ and $\pi_j=b$. We
build the permutation $\pi'$ corresponding to $\pi$ by setting
$\pi'_1=i$, $\pi'_n=j$, and $\pi_2'\ldots\pi_{n-1}'$ is obtained
from $\pi$ by first removing $a$ and $b$, and then, in what is left,
by replacing $i$ by $a$ and $j$ by $b$. For example, assuming that 2
and 4 are the determined elements, to
$13\underline{4}5\underline{2}6$ there corresponds
$\underline{5}1426\underline{3}$ which is a uniquely 5-determined
6-permutation.

Another application of Theorem~\ref{cr03} is finding lower and upper
bounds for the number $A_{k,n}$ of uniquely $k$-determined
$n$-permutations.

\begin{theorem}\label{number} We have $2((k-1)!)^{\lfloor n/k\rfloor}<A_{k,n}<2(2(k-1))^n$. \end{theorem}

\begin{proof}
According to Theorem~\ref{cr03}, we can estimate the number of
directed hamiltonian paths in $\mathcal{G}_{k,n}$ to get the
desired. This number is two times the number of (non-directed)
hamiltonian paths in $\mathcal{G}_{k,n}$, which is bounded from
above by $(2(k-1))^n$, since $2(k-1)$ is the maximum degree of
$\mathcal{G}_{k,n}$ (for $n\geq 2k-1$). So, $A_{k,n}<2(2(k-1))^n$.

To see that $A_{k,n}>2((k-1)!)^{\lfloor n/k\rfloor}$, consider
hamiltonian paths starting at node 1 and {\em not} going to any of
the nodes $i$, $i\geq k+1$ unless a path goes through {\em all} the
nodes $1,2,\ldots,k$. Going through all the first $k$ nodes can be
arranged in $(k-1)!$ different ways. After covering the first $k$
nodes we send the path under consideration to node $k+1$, which can
be done since we deal with $\mathcal{G}_{k,n}$. Then the path covers
{\em all}, but not any other, of the $k-1$ nodes $k+2,k+3,\ldots,
2k$ (this can be done in $(k-1)!$ ways) and comes to node $2k+1$,
etc. That is, we subdivide the nodes of $\mathcal{G}_{k,n}$ into
groups of $k$ nodes and go through all the nodes of a group before
proceeding with the nodes of the group to the right of it. The
number of such paths can be estimated from below by
$((k-1)!)^{\lceil n/k\rceil}$. Clearly, we get the desired result
after multiplying the last formula by 2 (any hamiltonian path can be
oriented in two ways).
\end{proof}

\subsection{Prohibitions giving unique
$k$-determinability}\label{proh}

The set of uniquely $k$-determined $n$-permutations can be described
by the language of prohibited patterns $\mathcal{L}_{k,n}$ as
follows. Using Theorem~\ref{cr01}, we can describe the set of
uniquely $k$-determined $n$-permutations by prohibiting patterns of
the forms $xX(x+1)$ and $(x+1)Xx$, where $X$ is a permutation on
$\{1,2,\ldots,|X|+2\}-\{x,x+1\}$ ($|X|$ is the number of elements in
$X$), the length of $X$ is at least $k-1$, and $1\leq x<n$. We
collect all such patterns in the set $\mathcal{L}_{k,n}$; also, let
$\mathcal{L}_k=\cup_{n\geq 0}\mathcal{L}_{k,n}$.

A prohibited pattern $X=aYb$ from $\mathcal{L}_k$, where $a$ and $b$
are some consecutive elements and $Y$ is a (possibly empty) word, is
called {\em irreducible} if the patterns of $Yb$ and $aY$ are not
prohibited, in other words, if the patterns of $Yb$ and $aY$ are
uniquely $k$-determined permutations. Without loss the generality,
we can assume that $\mathcal{L}_k$ consists only of irreducible
prohibited patterns.

\begin{theorem}\label{irr} Suppose $k$ is fixed. The number of (irreducible) prohibitions in $\mathcal{L}_k$ is
finite. Moreover, the longest prohibited patterns in $\mathcal{L}_k$
are of length $2k-1$.\end{theorem}

\begin{proof} Suppose that a pattern $P=xX(x+1)$ of length $2k$ or larger belongs to $\mathcal{L}_k$ (the case
$P=(x+1)Xx$ can be considered in the same way). Then obviously $X$
contains either $x-1$ or $x+2$ on the distance at least $k-1$ from
either $x$ or $x+1$. In any case, clearly we get either a prohibited
pattern $P'=yY(y+1)$ or $P'=(y+1)Yy$, which is a proper factor of
$P$. Contradiction with $P$ being irreducible.
\end{proof}

Theorem~\ref{irr} allows us to use the transfer matrix method to
find the number of uniquely $k$-determined permutations. Indeed, we
can consider the graph $\mathcal{P}_{2k-1}(\mathcal{L}_k)$, which is
the graph $\mathcal{P}_{2k-1}$ of patterns overlaps without nodes
containing prohibited patterns as factors. Then the number $A_{k,n}$
of uniquely $k$-determined $n$-permutation is equal to the number of
paths of length $n-2k+1$ in the graph, which can be found using the
transfer matrix method~\cite[Thm. 4.7.2]{Stan}\footnote{In fact, one
can use a smaller graph, namely $\mathcal{P}_{2k-2}(\mathcal{L}_k)$,
in which we mark arcs by corresponding permutations of length
$2k-1$; then we remove arcs containing prohibitions and use the
transfer matrix method. In this case, to an $n$-permutation there
corresponds a path of length $n-2k+2$. See Figure~\ref{il01} for
such a graph in the case $k=3$.}. In particular, the method makes
the following statement true.

\begin{theorem} The generating function $A_k(x)=\sum_{n\geq 0}A_{k,n}x^n$ for the number of uniquely $k$-determined
permutations is rational. \end{theorem}

A permutation is called {\em crucial} with respect to a given set of
prohibitions, if it does not contain any prohibitions, but adjoining
any element to the right of it leads to a permutation containing a
prohibition. In our case, an $n$-permutation is crucial if it is
uniquely $k$-determined, but adjoining any element to the right of
it, and thus creating an $(n+1)$-permutation, leads to a {\em
non}-uniquely $k$-determined permutation\footnote{As it is mentioned
in the introduction, crucial {\em words} are studied, for example,
in~\cite{EvdKit}. We define crucial permutations with respect to a
set of prohibited patterns in a similar way. However, as
Theorem~\ref{nonexist} shows, there are no crucial permutations with
respect to $\mathcal{L}_k$.}. If such a $\pi$ exists, then the path
in $\mathcal{P}_{2k-1}(\mathcal{L}_k)$ corresponding to $\pi$ ends
up in a sink. However, the following theorem shows that there are no
crucial permutations with respect to the set of prohibitions
$\mathcal{L}_k$, thus any path in
$\mathcal{P}_{2k-1}(\mathcal{L}_k)$ can always be continued.

\begin{theorem}\label{nonexist} There do not exist crucial permutations with respect to $\mathcal{L}_k$. \end{theorem}

\begin{proof} If $k=2$ then only the monotone permutations are uniquely $k$-determined, and we always can
extend to the right a decreasing permutation by the least element,
and the increasing permutation by the largest element.

Suppose $k\geq 3$ and let $Xx$ be an $n$-permutation avoiding
$\mathcal{L}_k$, that is, $Xx$ is uniquely $k$-determined. If $x=1$
then $Xx$ can be extended to the right by 1 without creating a
prohibition; if $x=n$ then $Xx$ can be extended to the right by
$n+1$ without creating a prohibition. Otherwise, due to
Theorem~\ref{cr01}, both $x-1$ and $x+1$ must be among the $k$
leftmost elements of $Xx$. In particular, at least one of them, say
$y$, is among the $k-1$ leftmost elements of $Xx$. If $y=x-1$, we
extend $Xx$ by $x$ (the ``old'' $x$ becomes $(x+1)$); if $y=x+1$, we
extend $Xx$ by $x+1$ (the ``old'' $x+1$ becomes $(x+2)$). In either
of the cases considered above, Theorem~\ref{cr01} guarantees that no
prohibitions will be created. So, $Xx$ can be extended to the right
to form a uniquely $k$-determined $(n+1)$-permutation, and thus $Xx$
is not a crucial $n$-permutation.
\end{proof}

\subsection{The case $k=3$}\label{k=3}

In this subsection we take a closer look to the graph
$\mathcal{P}_{4}(\mathcal{L}_3)$ whose paths give all uniquely
$3$-determined permutations (we read marked arcs of a path to form
the permutation corresponding to it). It turns out that
$\mathcal{P}_{4}(\mathcal{L}_3)$ has a nice structure (see
Figure~\ref{il01}).

Suppose $w'$ denotes the {\em complement} to an $n$-permutation
$w=w_1w_2\cdots w_n$. That is, $w'_i=n-w_i+1$ for $1\leq i\leq n$.
$\mathcal{P}_{4}(\mathcal{L}_3)$ has the following 12 nodes (those
are all uniquely $3$-determined $4$-permutations):\\

\begin{center}
\begin{tabular}{ll}  $a=1234$ & $a'=4321$ \\ $b=1324$ & $b'=4231$ \\ $c=1243$ & $c'=4312$ \\ $d=3421$ & $d'=2134$ \\
$e=1423$ & $e'=4132$ \\ $f=3241$ & $f'=2314$ \end{tabular}
\end{center}

In Figure~\ref{il01} we draw 20 arcs corresponding to the 20
uniquely $3$-determined $5$-permutations. Notice that
$\mathcal{P}_{4}(\mathcal{L}_3)$ is not strongly connected: for
example, there is no directed path from $c$ to $f$.

\setlength{\unitlength}{3mm}
\begin{figure}[h]
\begin{center}
\begin{picture}(20,22)
\put(0,0){\put(10,20){$a$} \put(10,14){$b$} \put(6,17){$d'$}
\put(14,17){$c$} \put(3,13){$f'$} \put(3,9){$f$} \put(6,5){$d$}
\put(14,5){$c'$} \put(10,2){$a'$} \put(10,8){$b'$} \put(17,13){$e$}
\put(17,9){$e'$}

\spline(9.9,20.6)(9.1,21.6)(10.3,22.8)(11.4,21.6)(10.8,20.6)
\spline(9.9,2)(9.1,1)(10.3,0)(11.4,1)(10.8,2)
\spline(2.8,10)(1.8,11.5)(2.8,13)\spline(4.2,10)(5.2,11.5)(4.2,13)
\spline(16.8,10)(15.8,11.5)(16.8,13)\spline(18.2,10)(19.2,11.5)(18.2,13)

\put(7.8,5.3){\vector(1,0){5.4}}\put(7.8,17.3){\vector(1,0){5.4}}

\put(4,8){\vector(1,-1){2}}\put(15,16){\vector(1,-1){2}}
\put(7.2,4.8){\vector(1,-1){2}}\put(11,20){\vector(1,-1){2}}
\put(7.2,18){\vector(1,1){2}}\put(11.2,2.5){\vector(1,1){2}}
\put(4,14.5){\vector(1,1){2}}\put(15,6.4){\vector(1,1){2}}
\put(9.5,14.7){\vector(-1,1){2}} \put(13,6){\vector(-1,1){2}}
\put(13.3,16.5){\vector(-1,-1){2}} \put(9.6,8){\vector(-1,-1){2}}

\shade\path(11.2,20.9)(10.8,20.6)(10.9,21)(11.2,20.9)\shade\path(10.8,1.6)(10.8,2)(11.2,1.7)(10.8,1.6)
\shade\path(2.4,10.3)(2.8,10)(2.8,10.5)(2.4,10.3)
\shade\path(4.2,12.5)(4.2,13)(4.7,12.7)(4.2,12.5)
\shade\path(16.4,10.3)(16.8,10)(16.8,10.5)(16.4,10.3)
\shade\path(18.2,12.5)(18.2,13)(18.7,12.7)(18.2,12.5) }

\end{picture}
\caption{Graph $\mathcal{P}_{4}(\mathcal{L}_3)$ (the case $k=3$).}
\label{il01}
\end{center}
\end{figure}
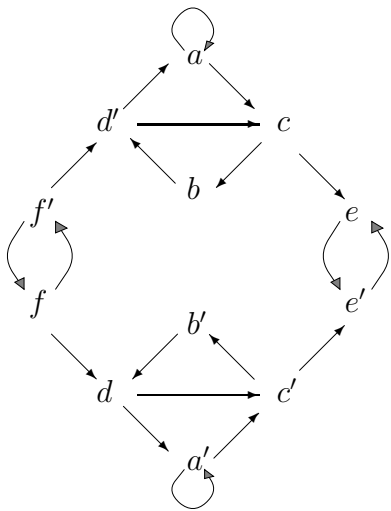

To find the generating function $A_3(x)=\sum_{n\geq 0}A_{3,n}x^n$
for the number of uniquely 3-determined permutations one can build a
12x12 matrix corresponding to $\mathcal{P}_{4}(\mathcal{L}_3)$ and
to proceed with the transfer matrix method. However, we do not do
that since, as it was mentioned in Subsection~\ref{dirpaths}, the
generating function for these numbers is
known~\cite[A003274]{Sloane}:
$$A_3(x)=\frac{1-2x+2x^2+x^3-x^5+x^6}{(1-x-x^3)(1-x)^2}.$$

\section{Open problems}\label{sec3}

It is clear that any $n$-permutation is uniquely $n$-determined,
whereas for $n\geq 2$ no $n$-permutation is uniquely $1$-determined.
Moreover, for any $n\geq 2$ there are exactly two uniquely
$2$-determined permutations, namely the monotone permutations. For a
permutation $\pi$, we define its {\em index} $IR(\pi)$ {\em of
reconstructibility} to be the minimal integer $k$ such that $\pi$ is
uniquely $k$-determined.

\begin{problem} Describe the distribution of $IR(\pi)$ among all $n$-permutations. \end{problem}

\begin{problem} Study the set of uniquely $k$-determined permutations in the case when a set of
nodes is removed from $\mathcal{P}_k$, that is, when some of
patterns of length $k$ are prohibited. \end{problem}

An $n$-permutation $\pi$ is $m$-$k$-{\em determined}, $m,k\geq 1$,
if there are exactly $m$ (different) $n$-permutations having the
same path in $\mathcal{P}_k$ as $\pi$ has. In particular, the
uniquely $k$-determined permutations correspond to the case $m=1$.

\begin{problem}\label{pr03} Find the number of $m$-$k$-determined $n$-permutations. \end{problem}

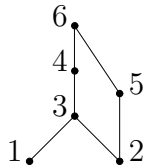
\begin{figure}[h]
\begin{center}
\begin{picture}(6,8)
\put(0,0){ \put(1,1){\p} \put(5,1){\p} \put(3,3){\p} \put(3,5){\p}
\put(3,7){\p} \put(5,4){\p} \put(0,1){1} \put(5.4,1){2} \put(2,3){3}
\put(2,5){4} \put(2,7){6} \put(5.4,4){5}
\path(1,1)(3,3)(5,1)\path(3,3)(3,5)(3,7)(5,4)(5,1)}
\end{picture}
\caption{The poset associated with the path $w=134265$ in
$\mathcal{P}_3$ ($k=3$).}\label{il03}
\end{center}
\end{figure}

Problem~\ref{pr03} is directly related to finding the number of {\em
linear extensions} of a {\em poset}. Indeed, to any path $w$ in
$\mathcal{P}_k$ there naturally corresponds a poset $\mathcal{W}$.
In particular, any factor of length $k$ in $w$ consists of
comparable to each other elements in $\mathcal{W}$. For example, if
$k=3$ and $w=134265$ then $\mathcal{W}$ is the poset in
Figure~\ref{il03}.

If all the elements are comparable to each other in $w$, then
$\mathcal{W}$ is a linear order and $w$ gives a uniquely
$k$-determined permutation. If $\mathcal{W}$ contains exactly one
pair of incomparable elements, then $w$ gives (two) 2-$k$-determined
permutations. In the example in Figure~\ref{il03}, there are 4 pairs
of incomparable elements, (1,2), (1,5), (3,5), and (4,5), and this
poset can be extended to a linear order in 7 different ways giving
(seven) 7-3-determined permutations.

\begin{problem} Which posets on $n$ elements appear while considering paths (of length $n-k+1$) in
$\mathcal{P}_k$? Give a classification of the posets (different from
the classification by the number of pairs of incomparable elements).
\end{problem}

\begin{problem} How many linear extensions can a poset (associated to a path in
$\mathcal{P}_k$) on $n$ elements with $t$ pairs of incomparable
elements have? \end{problem}

\begin{problem} Describe the structure of $\mathcal{L}_k$
(see Subsection~\ref{proh} for definitions) that consists of
irreducible prohibitions. Is there a nice way to generate
$\mathcal{L}_k$? How many elements does $\mathcal{L}_k$
have?\end{problem}


\begin{thebibliography}{11}
\bibitem{BurKit} A. Burstein and S. Kitaev: On unavoidable sets of word patterns, {\em SIAM J. on Discrete Math.} {\bf 19} (2005) 2,
371--381.
\bibitem{Chung} F. Chung, P. Diaconis and R. Graham: Universal
cycles for combinatorial structures, \emph{Discrete Math.} {\bf 110}
(1992), 43--60.
\bibitem{EhrKitPer} R. Ehrenborg, S. Kitaev and P. Perry: A Spectral Approach to Pattern-Avoiding Permutations, {\em 18th
International Conference on Formal Power Series \& Algebraic
Combinatorics}, the University of California, San Diego, USA, June
19--23 (2006).
\bibitem{ElizNoy} S. Elizalde and M. Noy: Consecutive subwords in permutations, {\em Advances in Applied Mathematics} {\bf 30} (2003), 110--125.
\bibitem{EvdKit} A. Evdokimov and S. Kitaev: Crucial
words and the complexity of some extremal problems for sets of
prohibited words, {\em Journal of Combinatorial Theory - Series A}
{\bf 105/2} (2004), 273--289.
\bibitem{GoulJack} I. P. Goulden and D. M. Jackson, {\em Combinatorial Enumeration}, A Wiley-Interscience Series in
Discrete Mathematics, John Wiley \& Sons Inc., New York, (1983).
\bibitem{Hur} G. Hurlbert: Universal Cycles: On Beyond de Bruijn,
PhD thesis, Department of Mathematics, Rutgers University, 1990.
\bibitem{Kit0} S. Kitaev: Multi-Avoidance of Generalised Patterns, {\em Discrete Math.} {\bf 260} (2003), 89--100.
\bibitem{Kit1} S. Kitaev: Partially ordered generalized patterns, {\em Discrete Math.} {\bf 298} (2005), 212--229.
\bibitem{Kit2} S. Kitaev: Counting independent sets on path-schemes, {\em Journal of Integer Sequences} {\bf 9},
no. 2 (2006), Article 06.2.2, 8pp.
\bibitem{Page} E. S. Page: Systematic generation of ordered sequences using recurrence relations,
{\em Computer Journal} {\bf 14} (1971), 150--153.
\bibitem{Sloane}
{\sc N. J. A.~Sloane and S.~Plouffe}, {\em The Encyclopedia of
Integer Sequences}, Academic Press, (1995)
\\ http://www.research.att.com/$\sim$njas/sequences/.
\bibitem{Stan} R. Stanley: {\rm Enumerative Combinatorics}, vol. I,
Cambridge Univ. Press, Cambridge, 1997.
\end{thebibliography}
\end{document}